\documentclass[12pt]{article}

      \usepackage{latexsym}
         \usepackage[reqno, namelimits, sumlimits]{amsmath}
         \usepackage{amssymb, amsfonts}
         \usepackage{amsthm}
 \newcommand{\beqn}{\begin{eqnarray}}
 \newcommand{\eeqn}{\end{eqnarray}}
 \newcommand{\be}{\begin{equation}}
 \newcommand{\ee}{\end{equation}}
 \newcommand{\ba}{\begin{array}{l}}
 \newcommand{\ea}{\end{array}}
 \newcommand{\pa}{\partial}
 \newcommand{\re}{\ref}
 \newcommand{\ci}{\cite}
 \newcommand{\la}{\label}

\newcommand{\ga}{\gamma}

\newcommand{\na}{\nabla}

 \newcommand{\De}{\Delta}

\newcommand{\M}{\mathcal{M}}

\renewcommand{\d}{{\mbox{div}}_g}
\newcommand{\D}{ {\Delta_g}}
\newcommand{\nag}{\nabla_g}
\newcommand{\nax}{\nabla_x}
\newcommand{\Dx}{\Delta_x}
\newcommand{\dx}{{\mbox{div}}_x}

\def\R{{\rm I\kern-.1567em R}}
\def\M{{\rm I\kern-.1567em M}}

\def\div {{\rm div\,}}

 \newtheorem{theorem}{Theorem}[section]
 
 \newtheorem{definition}[theorem]{Definition}
 
 \newtheorem{lemma}[theorem]{Lemma}
 
 \newtheorem{remark}[theorem]{Remark}

 \newtheorem{pro}[theorem]{Proposition}

\title{Global regularity of solutions of coupled Navier-Stokes equations and nonlinear Fokker Planck equations}

 \author{Peter Constantin\footnote{The University of Chicago} \and
  Gregory Seregin\footnote{Oxford University}}

\begin{document}
\maketitle

  \vspace{1cm}
 \noindent
 {{\bf Abstract } We provide a proof of global regularity of solutions of coupled Navier-Stokes equations and Fokker-Planck equations, in two spatial dimensions, in the absence of boundaries. The proof yields a priori estimates for the growth of spatial gradients.}

 \vspace {1cm}

\noindent {\bf 1991 Mathematical subject classification (Amer. Math.
Soc.)}: 35K, 35Q30, 82C31, 76A05.


\setcounter{equation}{0}
\section{Introduction  }
We consider a system
\be
\left \{
\ba
\pa_t u + u\cdot\nax u  - \nu\Dx u + \nax p = \dx \sigma, \\
\dx u = 0,\\
\pa _t f + u\cdot\nax f + \d( Wf) = \kappa (\D f + \d(f\nag U))
\ea
\right .
\la{sys0}
\ee
The functions $u= u(x,t)\in {\mathbb R}^2$, $p(x,t)\in \mathbb R$ and $f(x,m,t)$ are the unknown functions with $x\in{\mathbb T}^2$, $m\in M$ and $t\ge 0$ independent variables. $M$ is a compact connected Riemannian manifold without boundary of dimension $N$, with metric $g_{\alpha \beta}(m)$. The operations $\nag, \d, \D$ are covariant derivative of scalars, its adjoint and Laplace-Beltrami operator, respectively, i.e. in a local chart $m =(m^1, \dots m^N)$,
$$
\ba
\nag U = (\partial_{m^{\alpha}}U)_{\alpha = 1, \dots N}\\
\d W = \frac{1}{\sqrt{g}}\partial_{m^{\alpha}}\left (\sqrt{g}g^{\alpha \beta}W_{\beta}\right )\\
\D f = \d\nag f
\ea
$$
where, as customary, $g^{\alpha\beta}$ denotes the inverse of $g_{\alpha\beta}$,  $g$
its determinant and repeated indices are summed.
The second equation in (\re{sys0}) is the nonlinear Fokker-Planck equation.
Linear Fokker-Planck equations arise naturally as Kolmogorov forward equations for the probability density distributions associated with stochastic differential equations. Such linear equations with irregular coefficients were studied in (\cite{lebris}).
The cotangent field $W$ is given by
\be
W(x,m,t) = \nax u(x,t) : c(m) = \frac{\partial u_i(x,t)}{\partial x_j}c_{ji}(m),
\la{W}
\ee
where $c_{ji}(m) = (c_{ji; l}(m))_{l =1, \dots N}$ are smooth, time independent functions of $m$. The potential $U(x,m,t)$ is given by
\be
U(x,m,t) = \int_M k(m,m') f(x,m',t)dm' = U[f](x,m,t)
\la{U}
\ee
where $dm$ stands for the Riemannian volume element, $\sqrt{g}dm^1\dots dm^N$.
The interaction kernel $k$ is Lipschitz continuous in $M$ and is a given time independent function.  The kernel is a symmetric function $k(m,m') = k(m',m)$
and the operator $f\mapsto U[f]$ is bounded selfadjoint in $L^2(M)$.
The cotangent fields $c_{ji}$, the interaction kernel $k$,
the kinematic viscosity $\nu>0$ and the microscopic diffusivity $\kappa>0$ are all the parameters in the problem. The added stress $\sigma$ is given by them by
\be
\sigma_{ij}(x,t) = \int_M \left(c_{ij}\cdot \nag U(x,m,t) -\d c_{ij}\right)fdm.\la{sigma}
\ee
The system (\ref{sys0}) is a model of complex fluids (\cite{c}) that is a natural generalization of classical models of fluids with rod-like particles suspended in them (\cite{doied}). The constitutive relation (\ref{sigma}),  modulo multiples of the identity matrix, was  introduced (\cite{c-smo}) in order to have  the natural energy balance
\be
\frac{d}{dt}\int_{\mathbb R^d} \left \{\frac{1}{2} |u(x,t)|^2  + {\mathcal E}[f](x,t) \right\}dx + \int_{{\mathbb R}^d}\left \{\nu \left |\nax u\right|^2 + \kappa{\mathcal D}[f](x,t)\right\} dx = 0,
\la{bal}
\ee
where
\be
{\mathcal E}[f] = \int_M \left (f\log f + \frac{1}{2}fU[f]\right)dm
\la{fe}
\ee
and
\be
{\mathcal D}[f] = \int_M\left |\nag(U[f] +\log f)\right |^2fdm.
\la{dfe}
\ee
We note that the Fokker-Planck equations satisfies the weak maximum principle (if $u$ is a given smooth function) and therefore if $f(x,m, 0) \ge 0$, then $f(x,m,t) \ge 0$. Moreover the micoscopic density
$$
\rho_M(x,t) = \int_M f(x,m,t)dm
$$
obeys the pure transport equation
$$
\pa_t \rho_M  + u\cdot \nax \rho_M = 0,
$$
and therefore, the region occupied by microscopic corpora is material (carried by the flow). In particular, from the fact that $u$ is divergence-free it follows immediately that the density $\rho\in L^{\infty}(dt; L^1\cap L^{\infty}(dx))$, if the initial density $\rho_M(x,0)$ is bounded and integrable. The fluid density is taken to be identically $1$. We do not use the energy balance (\ref{bal}) in this paper, but the fundamental properties used for the proof of existence and regularity of solutions of (\ref{sys0}) originate from the same source as the energy balance, namely the structure of the equations and the constitutive equation (\ref{sigma}). In particular, we have, a priori,  $f\ge 0$, $f\in L^{\infty}(dx; L^1(M))$ and consequently $\sigma\in L^{\infty}$.

Global regularity for (\ref{sys0}) was first proved in (\cite{cm}). Independently and simultaneously, global regularity for a similar model was proved in
(\cite{fanghua}). That model is a version of the FENE model in which the physical gradient of velocity is replaced by its anti-symmetric part, and the
particles are restricted to the unit disk by a potential that is infinite at the unit circle. Both proofs suffer from the fact that they are  based on estimates with loss of regularity, and they are non-quantitative. In particular, there is no a priori bound on the growth rate for the spatial gradients.

In this paper we use results from (\cite{cs}) and the method of
(\cite{c-smo}) and (\cite{cftz}) to produce quantitative bounds. We use approximations
that respect the basic properties of (\ref{sys0}). We consider a standard mollifier
$$
J_{\delta}(u)(x) = \delta^{-2}\int_{\mathbb R^2} \phi\left(\frac{x-y}{\delta}\right )u(y,t)dy
$$
with $\phi\in C_0^{\infty}(\mathbb R^2)$ and approximate (\ref{sys0}) by
\be
\left \{
\ba
\pa_t u + u\cdot\nax u  - \nu\Dx u + \nax p = \dx J_{\delta}(\sigma), \\
\dx u = 0,\\
\pa _t f + J_{\delta}(u)\cdot\nax f + \d( J_{\delta}(W)f) = \kappa (\D f + \d(f\nag U))
\ea
\right .
\la{sys}
\ee
We obtain bounds independent of $\delta$. Existence of solutions of the nonlinear system (\ref{sys}) with $\delta>0$ can be obtained by an implicit iteration scheme, using linear equations in each step of the approximation:
\be
\left \{
\ba
\pa_t u^{(n+1)} + u^{(n))}\cdot\nax u^{(n+1)}  - \nu\Dx u^{(n+1)} + \nax p^{(n+1)} = \dx J_{\delta}\sigma^{(n+1)}, \\
\dx u^{(n+1)} = 0,\\
\pa _t f^{(n+1)} + J_{\delta}(u^{(n)})\cdot\nax f^{(n+1)} + \d(J_{\delta}(W^{(n)})f^{(n+1)}) = \\
\kappa (\D f^{(n+1)} + \d(f^{(n+1)}\nag U^{(n+1)}))
\ea
\right .
\la{sysn}
\ee
The existence of solutions of (\ref{sys0}) then follows from the existence of solutions of the approximate systems (\ref{sys}) and uniform bounds. The purpose of this article is to establish these bounds. In what follows all the bounds will be uniform in $\delta\ge 0$, and when $\delta =0$, then $J_\delta$ is taken to be the identity.

\begin{definition} Let $q>2$. We will say that $(u_0, f_0)$ are standard initial data if $\dx u_0 = 0$, $u_0\in W^{2,q}(\mathbb T^2)$, $f_0>0$,
$f_0(x,m)\in W^{1,q}(\mathbb T^2; L^2(M))$ and $\int_Mf_0(x,m)dm =1$.
\end{definition}
The Navier-Stokes equation can be written in the form
\be
\ba
\pa_t u -\nu\Dx u + \nax p= \dx\tau,
\\ \dx u = 0,
\ea
\la{nseq}
\ee
where
\be
\tau_{ij}(x,t) = J_{\delta}\sigma_{ij}(x,t) - u_i(x,t)u_j(x,t).\la{tau}
\ee
Taking the divergence of (\ref{nseq}) we can solve for the pressure
\be
p = -(-\Dx)^{-1}\pa_i\pa_j \tau_{ij}.\la{p}
\ee
using periodic boundary conditions. The operator $\tau\to p$ in (\ref{p}) is
bounded in $L^2$. The solutions of the Navier-Stokes equations discussed in this paper  have $\tau(\cdot, t)\in L^2$ and the pressure in the solution
of (\ref{sys}) is meant to be given by (\ref{p}).

\setcounter{equation}{0}
\section{Statements of the main result and lemmas}
The main result we prove in this paper is
\begin{theorem}\la{globounds}Let $q\ge 4$, $(u_0, f_0)$ be standard initial data and let
$T>0$ be arbitrary. Let $p>\frac{2q}{q-2}$, $\alpha>\frac{N}{2} +1$. There exists a constant $K$ depending only on the norms of the initial data, $T, \kappa, \nu, p, q, \alpha $, with $K$ bounded for bounded $T$, and a unique solution $(u,f)$ of (\ref{sys0}) with pressure $p$ given by (\ref{p}) and such that
\be
\|\nax\nax u\|_{L^p(0, T; L^q(\mathbb T^2))}\le K,
\la{ubound}
\ee
\be
\sup_{t\le T}\|\nax u(\cdot, t)\|_{L^{\infty}}\le K,
\la{naxt}
\ee
and
\be
\sup_{t\le T}\left \|f(\cdot, t)\right\|_{W^{1,q}(\mathbb T^2; H^{-\alpha}(M)))} \le K.
\la{fbound}
\ee
hold.
\end{theorem}
\begin{remark} The constant $K$ grows at most like a double exponential of $T$ multiplied by a first order polynomial in $T$.
\end{remark}
\noindent Let $H$ and $V$ be the completions of the set of all
divergence-free vector fields of $C^\infty(\mathbb T^2;\mathbb R^2)$
with vanishing mean value on the torus, with respect to the $L^2$
norm and the Dirichlet integral, respectively. The following result
was proved in \cite{cs}:
\begin{pro}\la{sercon} Let $u\in L^{\infty}(0,T;H)\cap L^2(0,T;V)$, $p\in L^2(0,T;L^2(\mathbb T^2))$ be a  solution of the initial value problem
\be\la{ns}\pa_tu+u\cdot\na u-\nu\De u+\na p=\div \sigma ,\qquad \div u=0,\ee
\be\la{id}u(\cdot,0)= u_0(\cdot)\in H,\ee
where $\sigma \in L^r(\mathbb T^2 \times (0,T);\mathbb M^{2\times 2})$ with $r\ge 4$.
Then, given $s>0$, there exists a constant $C_s$ depending only on $s$, $\nu$, the norm of $u_0$ in $H$, the norm of $\sigma$ in $L^r({\mathbb T}^2\times (0, T))$,  such that
\be\la{ub}\|u\|_{L^{\infty}(\mathbb T^2\times (s,T))}\le C_s.\ee
Moreover, the function $u$ is H\"older continuous in $\mathbb T^2\times [s,T]$ with exponent $\ga=1-\frac 4r$.\end{pro}
\begin{remark}\la{rns1} The existence  and uniqueness of a solution to the initial value problem
(\ref{ns}) and (\ref{id}) with above properties is well known, see \ci{La}.\end{remark}
The proof of the result is based on local iterative estimates for $L^4$ space-time integrals of the velocity, in the spirit of De Giorgi. The fact that $u, p$ solve the Navier-Stokes equation
(\ref{ns}), with $p$ given in (\ref{p}) with $\delta=0$, i.e.,
\be
p = R_iR_j(\sigma_{ij}-u_iu_j),\la{pu}
\ee
where $R_i = \partial_i(-\Dx)^{-{\frac{1}{2}}}$ are the Riesz transforms,
is used to relate the pressure to the velocity. The iteration relates integrals on smaller parabolic cubes to integrals on larger ones. For the iterative procedure to succeed, the modulus of absolute continuity of the map
$$
\Omega\subset \left\{\mathbb T^2 \times (0,T)\right\} \mapsto \int_{\Omega} |u(x,t)|^4dxdt,
$$
needs to be controlled apriori, to guarantee that such an integral is arbitrarily small, if the parabolic Lebesqgue measure of $\Omega$ is small enough.
The following result was used in (\cite{cs}) to control the modulus of absolute continuity.
\begin{pro}\label{sercon1}
Let $u\in L^{\infty}(0,T;H)\cap L^2(0,T;V)$  be a solution of the 2D Navier-Stokes equations (\ref{ns}) with initial data (\ref{id}), $u_0\in H\cap L^r({\mathbb T}^2)$ and  $\sigma \in L^r(\mathbb T^2 \times (0,T);\mathbb M^{2\times 2})$ with $r\ge 4$.
There exists a constant $K$ depending only on the norm $\|\sigma\|_{L^{r}(\mathbb T^2\times (0,T))}$,
$\nu, T$ and the norm of $u_0$ in $H\cap L^r(\mathbb T^2)$ such that
\be
\sup_{0 \le t\le T}\|u(\cdot,t)\|_{L^r(\mathbb T^2)} \le K.
\la{ur}
\ee
\end{pro}
\begin{pro}\label{localpro} Let $(u_0,f_0)$ be standard initial data. There exists $T_0 >0$ and a constant $K$, depending only on the initial data and the parameters $\nu, \kappa, p, q, \alpha$, where $q>2, p>\frac{2q}{q-2}, \alpha > 1 + \frac{N}{2}$,  such that a unique solution of (\ref{sys0}) exists on the time interval $[0, T_0]$, with pressure $p$ given by (\ref{p}), satisfying
\be
\|u\|_{L^p(0, T_0; L^q(dx))} + \sup_{0\le t\le T_0}[\|u(\cdot, t)\|_{W^{1,q}(\mathbb T^2)} + \|f(\cdot, t)\|_{W^{1,q}(\mathbb T^2; H^{-\alpha}(M))}] \le K.
\la{loc}
\ee
\end{pro}
\begin{pro}\la{linftypro} Let $(u_0, f_0)$ be standard initial data and let $T>0$ be arbitrary. Let $q\ge 2$, $p>\frac{2q}{q-2}$, $\alpha>\frac{N}{2} +1$. There exists a constant $K$ depending only on the initial data, $T, \kappa, \nu, p,q,\alpha$, bounded for bounded $T$, such that, if $(u,f)$ is a solution of (\ref{sys}) or a solution of (\ref{sys0}) with
$$
\ba
u\in L^p(0,T; W^{2,q}(\mathbb T^2)), \; {\mbox{pressure}}\;p\; {\mbox{given by}}\; (\ref{p}), \\ \;\;{\mbox{and with}}\;\;
f\in L^p(0,T; W^{1,q}(\mathbb T^2; H^{-\alpha}(M))),
\ea
$$
then
\be
\sup_{0\le t\le T}\left [\|u (\cdot,t )\|_{L^{\infty}(\mathbb T^2)}^2 +
\|\sigma (\cdot, t)\|_{L^{\infty}(\mathbb T^2)}\right ] \le K
\la{K}
\ee
holds.
\end{pro}
\begin{lemma}\la{gradlemma} Let $(u_0, f_0)$ be standard initial data and let $T>0$ be arbitrary. Let  $q\ge 4$, $p>\frac{2q}{q-2}$, $\alpha>\frac{N}{2} +1$. There exists a constant $K$ depending only on the initial data, $T, \kappa, \nu, p,q,\alpha $, bounded for bounded $T$, such that, if $(u,f)$ is a solution of (\ref{sys}) or  a solution of (\ref{sys0}) with $u\in L^p(0,T; W^{2,q}(\mathbb T^2))$, pressure $p$ given by (\ref{p}) and
$f\in L^p(0,T; W^{1,q}(\mathbb T^2; H^{-\alpha}(M)))$, then
\be
\sup_{0\le t\le T}\|\nax u(\cdot,t)\|_{L^{\infty}(\mathbb T^2)}\le K\log( 2 +
\left \|f\right\|_{L^p(0,T; W^{1,q}(\mathbb T^2; H^{-\alpha}(M)))}).
\la{basicg}
\ee
\end{lemma}

\begin{theorem} {\la{main}} Let $T>0$, and $u_0, f_0$ be arbitrary standard initial data. Let $k\ge 1$ $q\ge 4$ and assume that $u_0\in W^{k+1,q}(\mathbb T^2)$ and $f_0\in W^{k,q}({\mathbb T^2}; L^q(M))$.
Then, for any $p>\frac{2q}{q-2}$, there exist constants $K$
depending only on $k, q, p$, $\nu,\kappa, T$ and the norms of the initial data, such that the solution of (\ref{sys0}) on $[0,T]$ with pressure $p$ given by (\ref{p}) sstisfies
\be
\sup_{0\le T}\|u(\cdot, t)\|_{W^{k,q}(\mathbb T^2)} + \|u\|_{L^p(0,T; W^{k+1,q}(\mathbb T^2))} \le K
\la{basick}
\ee
and
\be
\sup_{t\le T}\|f\|_{W^{k,q}(\mathbb T^2; L^q(M)))}\le K.
\la{fk}
\ee

\end{theorem}

\section{Proof of Proposition \ref{sercon1}}
We multiply (\ref{ns}) by $u|u|^{r-2}$ and integrate in space and integrate by parts. We obtain
$$
\ba
\frac{d}{rdt}\int_{\mathbb T^2}|u(x,t)|^rdx  + \nu\int_{\mathbb T^2} |\nax u(x,t)|^2|u(x,t)|^{r-2}dx \le\\
C\int_{\mathbb T^2}\left[|\sigma(x,t)| + |p(x,t)|\right]|\nax u(x,t)||u(x,t)|^{r-2}dx.
\ea
$$
Writing $|\nax u||u|^{r-2} = |\nax u||u|^{\frac{r-2}{2}}|u|^{\frac{r-2}{2}}$,
using a H\"{o}lder inequality with exponents $r, 2,  \frac{2r}{r-2}$, and then the Schwartz inequality and the viscous term, we deduce
$$
\frac{d}{dt}\left\|u(\cdot,t)\right\|^r \le C\left [\|p(\cdot,t)\|_{L^r}^2 + \|\sigma (\cdot, t) \|_{L^r}^2\right ]\|u\|_{L^r}^{r-2}
$$
Dividing by $\|u\|_{L^r}^{r-2}$ and using the boundedness of Riesz transforms we obtain
$$
\frac{d}{dt}\|u\|_{L^r}^2 \le C\left [\|u\|_{L^{2r}}^4 + \|\sigma\|_{L^r}^2\right ]
$$
Now we use the inequality
\be
\|u\|^2_{L^{2r}({\mathbb T^2})} \le C\left [\|u\|_{L^2(\mathbb T^2)} +\|\nax u\|_{L^2(\mathbb T^2)}\right ]\|u\|_{L^r(\mathbb T^2)}
\la{interp}
\ee
to deduce that
$$
\frac{d}{dt}\|u\|_{L^r}^2 \le C\|\sigma\|_{L^r}^2 + C\left [\|u\|_{L^2}^2 + \|\nabla u\|^2_{L^2}\right]\|u\|_{L^r}^2.
$$
Because
$$
\int_0^T \left [\|u\|^2_{L^2} + \|\nabla u\|^2_{L^2}\right ]dt \le C
$$
the bound (\ref{ur}) follows from Gronwall's inequality. We present below a sketch of the proof of (\ref{interp}).
We claim first that for any $r\ge 2$ there exists a constant $C_r$, such that
\be
\|f\|_{L^{2r}(\mathbb R^2)}^2 \le C_r\|f\|_{L^{r}(\mathbb R^2)}\|\nabla f\|_{L^2(\mathbb R^2)}
\la{fr}
\ee
holds for all $f\in L^r(\mathbb R^2)$ with $\nabla f \in L^2(\mathbb R^2)$.
This is a generalization of the well-known Ladyzhenskaya inequality (\cite{La}) corresponding to $r=2$. An elemenatry proof of (\ref{fr}) was given in (\cite{cs}). We give, for the sake of completeness, a generalization and proof in the Appendix.
The inequality (\ref{interp}) follows by considering the function $u$ as the restriction to $[0,2\pi]^2$ of a periodic function $U$ defined in the whole space $\mathbb R^2$,  and taking a smooth compactly supported function $\phi$ that is identically $1$ on an open neighbourhood of $[0,2\pi]^2$. The inequality (\ref{fr}) holds for $f = \phi U$, and in view of the fact that
$\|\phi U\|_{L^r(\mathbb R^2)} \le C\|u\|_{L^r(\mathbb T^2)}$ and similar inequalities, we have
$$
\ba
\|u\|^2_{L^{2r}(\mathbb T^2)} \le \|f\|^2_{L^{2r}(\mathbb R^2)} \le\\
C\|\nabla (\phi U)\|_{L^2(\mathbb R^2)}\|\phi U\|_{L^r(\mathbb R^2)}\le\\
 C\left [\|\phi\nabla U\|_{L^2(\mathbb R^2)} + \|U\nabla \phi\|_{L^2(\mathbb R^2)}\right ]\|u\|_{L^r(\mathbb T^2)}\le\\
C\left [\|\nabla u\|_{L^2(\mathbb T^2)} + \|u\|_{L^2(\mathbb T^2)}\right]\|u\|_{L^r(\mathbb T^2)}.
\ea
$$
\section{Proof of Proposition \ref{localpro} }
In this section we will denote by $C$ constants that may depend on $\nu, \kappa, T, p, q$ and $\alpha$ and are locally bounded in $T>0$. We will denote by $K$ constants that may depend in addition on standard initial data, and are locally bounded in $T$ and the norms of standard initial data. All the constants are independent of  $\delta\ge 0$.

We consider the vorticity, $\omega (x,t) = \nabla^{\perp}\cdot u(x,t) =
\frac{\partial u^2(x,t)}{\partial x^1} - \frac{\partial u^1(x,t)}{\partial x^2}$
and, taking the curl of the Navier-Stokes equation, we obtain the vorticity equation
\be
\pa_t \omega + u\cdot\nax \omega  -\nu\Dx \omega = \nax^{\perp}\cdot\dx J_{\delta}\sigma.
\la{vorteq}
\ee
We multiply this by $|\omega|^{q-2}\omega$ and integrate in space:
$$
\ba
\frac{d}{qdt}\int_{\mathbb T^2}|\omega|^q dx + \nu(q-1) \int_{\mathbb T^2}\left |\nax \omega\right|^2 |\omega|^{q-2} dx\\
 \le (q-1)\int _{\mathbb T^2}|\dx J_{\delta}\sigma| |\omega|^{\frac{q-2}{2}}\left [|\nax\omega||\omega|^{\frac{q-2}{2}}\right]dx.
\ea
$$
Using the H\"{o}lder inequality with exponents $q, \frac{2q}{q-2}, 2$
and hiding the term involving the gradient of $\omega$ in the viscous term, we obtain
\be
\ba\frac{d}{qdt}\int_{\mathbb T^2}|\omega|^q dx + \frac{\nu(q-1)}{2} \int_{\mathbb T^2}\left |\nax \omega\right|^2 |\omega|^{q-2} dx\\
 \le \frac{q-1}{2\nu}\|\dx J_{\delta}\sigma\|_{L^q}^2\|\omega\|_{L^q}^{q-2}.
\ea
\la{omegineq}
\ee
Integrating and using the well-known fact that the $L^q$ norms of vorticity bound from above the $L^q$ norms of the full gradient of velocity (modulo multiplicative constants),  we obtain
\be
\|\nax u(\cdot ,t)\|_{L^q}^2 \le K + C\int_0^t\|\dx J_{\delta}\sigma(\cdot, s)\|^2_{L^q}ds.
\la{naxsigq}
\ee
The forces applied by the particles are obtained after $f$ is integrated along with smooth coefficients on $M$ in order to produce $\sigma$ (\ref{sigma}).
 Therefore, only very weak regularity of $f$ with respect to the microscopic variables $m$ is sufficient to control $\sigma$.
We take advantage of this fact in order to control spatial gradients of $f$
in terms of $\nax u\in L^1(L^{\infty})$.  We consider the $L^2(M)$ selfadjoint pseudodifferential operator
\be
R = \left (-\D +{\mathbf I}\right)^{-\frac{\alpha}{2}}\la{R}
\ee
with $\alpha> \frac{N}{2} + 1$.
We differentiate the Fokker-Planck equation
\be
\pa_t f + J_{\delta}u\cdot\nax f +\d(J_{\delta}(W)f) = \kappa\d(\nag(\log f + U[f]))\la{feq}
\ee
in (\ref{sys}) with respect to $x$, apply $R$, multiply by
$R\nax f$ and integrate on M. Let us denote by
\be
N(x,t)^2 = \int\limits_{M}\left | R\nax f (x,m,t)\right|^2dm \la{Nx}
\ee
the square of the $L^2$ norm of $R\nax f$ on M. Note that
\be
|\nax J_{\delta}\sigma(x,t)| \le C N(x,t)\la{naxsign}
\ee
holds in view of the definition (\ref{sigma}). We obtain
\be
\frac{1}{2}\left(\partial_t + J_{\delta}u\cdot \nax \right )N^2 \le C(|J_{\delta} \nax u| +\kappa)N^2 + C|J_{\delta}\nax\nax u|N
\la{Nineq}
\ee
pointwise in $(x,t)$ with an absolute constant C. The proof of this fact
appeared in several places (\cite{c}, \cite{c-smo}, \cite{cftz}) and will not be reproduced here.
Now we multiply  (\ref{Nineq}) by $N^{q-2}$, integrate $dx$,
multiply by $\|N\|_{L^q(dx)}^{p-q}$ and use  H\"{o}lder inequalities in both space and  time:
\be
\ba
\frac{d}{pdt}\|N(\cdot,t)\|^p_{L^q(dx)} \le \\
\le C \|J_{\delta}\nax\nax u(\cdot, t)\|^p_{L^q(dx)} + C\left(\|J_{\delta}\nax u(\cdot, t)\|_{L^{\infty}(dx)} + 1\right )\|N(\cdot, t)\|^p_{L^q(dx)}.
\ea
\la{nlpq}
\ee
In order to proceed we need to use the representation of the gradient of the solution, from the Navier-Stokes equation (\ref{nseq}):
\be
\nax u (x,t) = e^{\nu t\Delta}\nax u_0  - \int\limits_0^te^{\nu(t-s)\Delta} \Delta{\mathbb H}\tau(\cdot, s)ds
\la{naxhst}
\ee
with
\be
({\mathbb H}\tau)_{ij} = R_j\left (\delta_{il} + R_iR_l\right)R_k\tau_{lk}
\la{hsigma}
\ee
and $R_j =\pa_j(-\Dx)^{-\frac{1}{2}}$ are Riesz transforms.
This formula is obtained by differentiation of the integral representation of
the solution of the Navier-Stokes equation (\ref{nseq}) and use of (\ref{p}).
We will use the fact that the linear operator
$$
h(t) \mapsto {\mathcal T}h = \int_0^t e^{\nu(t-s)\Delta}\Delta {\mathbb H}h(s)ds
$$
is bounded in $L^p(dt; L^q(dx))$ for $1<p, q<\infty$ (see, for example
(\cite{lemarie})). We start by estimating $\nabla u$ from (\ref{naxhst}) using the smoothness of the kernel of the heat equation  which results in the bound (see also (\cite{c-smo}))
$$
\|e^{\nu(t-s)\Delta}\Delta{\mathbb H}\tau(s)\|_{L^{\infty}(dx)}\le C(t-s)^{-1}\|\tau (s)\|_{L^{\infty}(dx)},
$$
and obtain, for any $0<l<t$,
$$
\ba
\|\nax u(\cdot, t)\|_{L^{\infty}(dx)} \le K + C\|\tau\|_{L^{\infty}(0;t)}\int\limits_0^{t-l}(t-s)^{-1}ds\\
+ \int\limits_{t-l}^t\|e^{\nu(t-s)\Delta}\dx{\mathbb H}\nax\tau(s)\|_{L^{\infty}(dx)}ds.
\ea
$$
We measure $\nax\tau$ in  $L^p(0,t; L^q(dx))$ for $p>\frac{2q}{q-2}$, $q>2$.
Then we obtain, using properties of the heat kernel,
$$
\ba
\|\nax u (\cdot, t)\|_{L^{\infty}}\le K + C\|\tau\|_{L^{\infty}(0,t)}\log\left (\frac{t}{l}\right ) \\
+ C\int\limits_{t-l}^t(t-s)^{-\frac{q+2}{2q}}\|\nax\tau(\cdot,s)\|_{L^q(dx)}ds
\ea
$$
and thus, by the H\"{o}lder inequality in time with $p, p^*$ we obtain that the  last term is bounded by
$$
l^{(\frac{1}{p^*}-\frac{q+2}{2q})}Y_{pq}(t)
$$
with
\be
Y_{pq}(t) = \left (\int_0^t\|\nax\tau (\cdot, s)\|_{L^q(dx)}^pds\right)^{\frac{1}{p}}.
\la{Y}
\ee
The choice of $p$ was so that the power of $l$ is positive. Then,  choosing
$l$ in terms of $Y_{pq}$ we get:
\be
\|\nax u(\cdot,t)\|_{L^{\infty}(dx)} \le K + C\|\tau\|_{L^{\infty}(0,t)}\log \left( 2 + Y_{pq}(t)\right)
\la{naxuhstbound}
\ee
for $0\le t\le T$.
From (\ref{naxhst}) we deduce
$$
\nax\nax u = e^{\nu t\Delta}\nax\nax u_0  + \int\limits_0^te^{\nu(t-s)\Delta}\Delta ({\mathbb H}\nax\tau(s))ds
$$
and therefore, using the boundedness of ${\mathcal T}$ in $L^p(dt;L^q(dx))$, we deduce
\be
\|\nax\nax u\|_{L^p((0,t), L^q({\mathbb T}^2))} \le K + CY_{pq}(t).
\la{naxnaxbound}
\ee
Let us denote
\be
Z_{pq}(t) = \left (\int_0^t\|N(\cdot, s)\|_{L^q}^pds\right)^{\frac{1}{p}}\la{Z}
\ee
Then, using (\ref{naxuhstbound}) and (\ref{naxnaxbound}) in (\ref{nlpq}), integrating in time and using the fact that the right-hand side of (\ref{naxuhstbound}) is non-decreasing in time, we have
\be
\frac{d}{dt}Z_{pq}(t)^p\le K + CY_{pq}(t)^p + C\|\tau\|_{L^{\infty}(0,t)}\log(2 + Y_{pq}(t))Z_{pq}(t)^p
\la{zbound}
\ee
Now, from the definition of $Y_{pq}$, the Sobolev embedding $W^{1,q}(\mathbb T^2) \subset L^{\infty}(\mathbb T^2)$, (\ref{naxsigq}) and (\ref{naxsign}), it follows that
\be
Y_{pq}(t) \le C(Z_{pq}(t) + Z_{pq}^2(t))\la{YZ}.
\ee
Also, using (\ref{sigma}), the Sobolev embedding referred to above and (\ref{naxsigq}), it follows that
\be
\|\tau \|_{L^{\infty}} \le K + CZ_{pq}^2(t),
\la{taub}
\ee
and thus, using (\ref{YZ}) and (\ref{taub}) in (\ref{zbound}) we obtain an ordinary differential inequality for $Z_{pq}^p$ which shows that there exists $T_0$ and $K$ depending on the norms of standard initial data such that
$$
Z_{pq}(t) \le K, \quad {\mbox{for}}\;\; t\le T_0.
$$
Using (\ref{YZ}) we deduce that $Y_{pq}(t)$ is bounded a priori in terms of initial data, and using (\ref{naxuhstbound}), (\ref{naxnaxbound}) and  (\ref{taub}),  we obtain  a priori bounds for $u$. The proof of local existence follows by passing to the limit $\delta\to 0$.
\section{Proof of Proposition {\ref{linftypro}}}
In this section we have a solution of (\ref{sys})
or (\ref{sys0}) with standard initial data, with $u\in L^p(0,T; W^{2,q}(\mathbb T^2))$, $f \in L^p(0,T; W^{1,q}(\mathbb T^2; H^{-\alpha}(M)))$, $p$ given by (\ref{p}), and we need to find $K$ depending only on the initial data
and parameters $\nu, \kappa, q$ but not on the solution (i.e, bounded a priori in terms of the initial data and parameters) such that (\ref{K}) holds.

Note first of all that $\|J_{\delta}\sigma\|_{L^{\infty}(\mathbb T^2; \mathbb M^{2\times 2})} \le K$ follows from
(\ref{sigma}) with a uniform constant, depending only on the coefficients $c_{ij}$ and $k$ of the nonlinear Fokker Planck equation and on the initial density of particles, which in this paper we took to be 1.

It remains to bound $\|u\|_{L^{\infty}}^2$. We may apply the local existence result, Proposition \ref{localpro}.
Let $T_0>0$ be a local existence time depending on the standard data and parameters, guaranteed to exist by Proposition \ref{localpro}. If $T\le T_0$ there is nothing left to prove. If $T>T_0$, we take $s= \frac{T_0}{2}$ in Proposition \ref{sercon}. The assumptions of Proposition \ref{sercon} are satisfied: Indeed, the right-hand side $J_{\delta}\sigma$ is in $L^r$ with arbitrary $r$. In addition, the solution $u$ is sufficiently regular to justify the standard energy inequality
$$
\ba
\frac{1}{2}\|u(\cdot t)\|_{L^2}^2 + \frac{\nu}{2}\int_0^T\int_{\mathbb T^2}\left |\nax u(x,t)\right |^2dxdt \le
\\
\frac{1}{2}\|u_0\|^2_{L^2} + \frac{1}{2\nu}\int_0^T\int_{\mathbb T^2}|J_{\delta}\sigma (x,t)|^2dxdt
\ea
$$
and therefore, by Sobolev embedding $u\in L^4((0,T)\times {\mathbb T^2})$. In view of $p = R_iR_j(\tau_{ij})$ from (\ref{p}), and the boundedness of Riesz operators in $L^2$, the assumptions of Proposition {\ref{sercon}} are verified. Because
$s=\frac{T_0}{2}$ is fixed by the initial data, the constant $C_s$ in Proposition (\ref{sercon}) is bounded uniformly in terms of the norms of the standard initial data and $T$.

\section{Proof of Lemma {\ref{gradlemma}}, and Theorems {\ref{globounds}} and {\ref{main}}}
For the proof of Lemma \ref{gradlemma} we use (\ref{naxuhstbound}). In view of
Proposition \ref{linftypro} and the definition (\ref{tau}) we have
\be
\sup_{0\le t\le T}\|\tau (\cdot,t)\|_{L^{\infty}} \le K\la{tauk}
\ee
with $K$ depending only on initial data, parameters, an bounded locally in $T>0$.
We can use the bound (\ref{YZ}) together with (\ref{tauk}) in (\ref{naxuhstbound}) to obtain
\be
\|\nax u (\cdot, t)\|_{L^{\infty}(\mathbb T^2)} \le K\log(2 +Z_{pq}(t))
\la{naxuz}
\ee
and, consequently, (\ref{basicg}) and Lemma \ref{gradlemma} are proved. In order to prove Theorem \ref{globounds} we use the bound (\ref{K}) of Proposition \ref{linftypro}, the chain rule, the bound (\ref{naxsigq}) and (\ref{naxsign}) to obtain
\be
\|\nax \left(u\otimes u\right)\|_{L^p(0,t; L^q(\mathbb T^2))} \le K(1+ Z_{pq}(t)),
\la{naxuu}
 \ee
for $t\le T$. On the other hand, from (\ref{naxsign}) we have
$$
\|\nax \sigma\|_{L^p(0,t; L^q(\mathbb T^2))} \le CZ_{pq}(t),
$$
so we obtain, in view of the definition (\ref{Y}),
\be
Y_{pq}(t) \le K(1+Z_{pq}(t))\la{YZlin}
\ee
for $t\le T$. In view of (\ref{naxnaxbound}) and (\ref{YZlin}) above, we have
\be
\|\nax\nax u\|_{L^p(0,t; L^q(\mathbb T^2))} \le K(1 + Z_{pq}(t))\la{naxnaxlin}\ee
for $0\le t\le T$, as well. Using (\ref{tauk}) and (\ref{YZlin}) in (\ref{zbound}) we deduce that $z(t)= Z_{pq}^p(t)$ obeys an ordinary differential inequality
$$
\frac{d}{dt}z \le K(2+z)\log(2+z)
$$
with $z(0) = 0$, $z\ge 0$ and $K(t)$ locally bounded on $[0, \infty)$. This implies an apriori bound for $2+z(t)$,
$$
2+z(t)\le \exp{\left[(\log 2)\exp\left (\int_0^tK(s)ds\right)\right]}.
$$
Thus
\be
Z_{pq}(T)\le K.
\la{ZT}
\ee
It follows from (\ref{ZT}) and (\ref{naxnaxlin})
that $\|\nax\nax u\|_{L^p(0, T; L^q(\mathbb T^2))}\le K$, which proves (\ref{ubound}). Using (\ref{naxuz}) and (\ref{ZT}), we obtain (\ref{naxt}).
Finally, using (\ref{ubound}) and (\ref{naxt}) in (\ref{nlpq}), we obtain (\ref{fbound}). The proof of Theorem \ref{globounds}
is complete.

The proof of Theorem {\ref{main}} is done by induction. For $k=1$ , the inequality (\ref{basick}) follows from (\ref{ubound}) and (\ref{naxt}). The inequality
(\ref{fk}) follows from the Fokker-Planck equation (\ref{feq}). Indeed, differentiating with respect to $m$, using the fact that the coefficient $\kappa>0$ of $\D$ is nonzero and the bound (\ref{naxt}), we obtain first that $\nag f\in L^{\infty}(0,T;L^{\infty}(dx; L^q(M)))$.  We also obtain an apriori $f\in L^{\infty}$ bound, directly from (\ref{feq}) using (\ref{naxt}). Then, we differentiate (\ref{feq}) in $x$, multiply by $\nax f |\nax f|^{q-2}$ and integrate $dm$. We obtain, after integrations by parts,
$$
\ba
\left (\pa _t + J_{\delta}u(x,t)\cdot\nax\right) \|f(x, \cdot, t)\|_{L^q(M)}^q \le\\
- \kappa\int_M |\nag\nax f(x, m,t)|^2 |\nax f(x,m,t)|^{q-2}dm + \\
C|\nax J_{\delta}u(x,t)|\|f(x,\cdot, t) \|_{L^q(M)}^q + \\
|J_{\delta} \nax\nax u(x,t)|\|\d(cf)\|_{L^q(M)}\|f(x,\cdot, t)\|_{L^q(M)}^{\frac{q-1}{q}} +\\
C\int_M|\nag U||\nax f|^{q-1}|\nag\nax f|dm + \\
C\int_M|f||\nag \nax U||\nax f|^{q-2}|\nag\nax f|dm.
\ea
$$
The last term is bounded using the dissipative term involving $\kappa$,
$$
\ba
\int_M|f||\nag \nax U||\nax f|^{q-2}|\nag\nax f|dm \le\\
\frac{\kappa}{2}\int_M |\nag\nax f(x, m,t)|^2 |\nax f(x,m,t)|^{q-2}dm\\
+C\|f(x, \cdot,t)\|_{L^{\infty}(M)}^2\|\nax f(x,\cdot,t)\|^q_{L^q(M)}.
\ea
$$
We do have $\|f(x, \cdot,t)\|_{L^{\infty}(M)}^2 \le K$ because of (\ref{naxt}).
We bound the penultimate term similarly:
$$
\ba                                                                            \int_M|\nag U||\nax f|^{q-1}|\nag\nax f|dm \le\\
\frac{\kappa}{2}\int_M |\nag\nax f(x, m,t)|^2 |\nax f(x,m,t)|^{q-2}dm\\        +C\|\nag U(x, \cdot,t)\|_{L^{\infty}(M)}^2\|\nax f(x,\cdot,t)\|^q_{L^q(M)}.    \ea                                                                             $$
The term $\|\nag U\|_{L^{\infty}(M)}\le K$ because $f (x,\cdot, t)$ is bounded in $L^1(M)$ annd the kernel $k(m,m')$ is Lipschitz. Integrating in space
and using (\ref{naxt}) we obtain
$$
\frac{d}{dt}\|f\|_{L^q(dx; L^q(M))}\le K\|\nax\nax u\|_{L^q}
$$
the bound (\ref{fk}) for $k=1$. For the induction step, we differentiate  (\ref{naxhst}) $k$ times. Using the classical calculus inequality
$$
\|u\otimes u\|_{W^{k,q}} \le C\|u\|_{L^{\infty}}\|u\|_{W^{k,q}}
$$
and the induction hypothesis (\ref{basick}) and (\ref{fk}), we obtain
(\ref{basick}) for $k+1$. In order to obtain (\ref{fk}) for $k+1$ we
differentiate (\ref{feq}) $k+1$ times, use the fact that (\ref{basick}) is true for $k+1$, and employ arguments similar to those shown for $k=1$. We omit further details.
\section{Appendix: Generalized Ladyzhenskaya inequalities}
The inequality (\ref{fr}) is the particular case $n=2$ of the inequality
\be
\|f\|_{L^{2r}(\mathbb R^n)}^2 \le C\|f\|_{L^r(\mathbb R^n)}\|\nabla f\|_{B^{0, n}_2(\mathbb R^n)}
\la{fro}
\ee
valid for all $r\ge \frac{n}{2}$. The norm in right the hand side of (\ref{fro}) is the Besov space norm,
$$
\|f\|_{B^{s,p}_q(\mathbb R^n)} = \left [\sum_{j=-\infty}^\infty \lambda_j^{qs}\|\Delta_j f\|_{L^p(\mathbb R^n)}^q\right ]^{\frac{1}{q}}
$$
defined in terms of the  Littlewood-Paley decomposition (\cite{bl}, \cite{lemarie})
$$
f = \sum_{j\in {\mathbb Z}}\Delta_j f
$$
into functions whose Fourier transforms are supported in dyadic shells
of order $\lambda_j = 2^j$, i.e. $\sup\widehat{\Delta_j f}\subset A_j$,
$A_j=\{\xi\left | |\xi|\in [2^{j-1}, 2^{j+1}]\right .\}$.
As it is well known $B^{0, 2r}_2(\mathbb R^n) \subset L^{2r}(\mathbb R^n)$, i.e.,
$$
\|f\|_{L^{2r}(\mathbb R^n)}^2 \le C\sum_{j\in\mathbb Z}\|\Delta_j f\|_{L^{2r}(\mathbb R^n)}^2.
$$
We split the sum in two parts, for $j\le M$ and for $j\ge M$.
When $j\ge M$ we use  the Bernstein inequality
$$
\|\Delta_j f\|_{L^{2r}({\mathbb R^n})}^2\le C\lambda_j^{-\frac{n}{r}}\|\nabla\Delta_j f\|^2_{L^n(\mathbb R^n)}
$$
and thus
$$
\ba
\sum_{j\ge M}\|\Delta_j f\|^2_{L^{2r}(\mathbb R^n)} \le C\sum_{j\ge M}\lambda_j^{-\frac{n}{r}}\|\Delta_j\nabla f\|^2_{L^n(\mathbb R^n)}\le \\
C\lambda_M^{-\frac{n}{r}}\|\nabla f\|_{B^{0,n}_2(\mathbb R^n)}^2.
\ea
$$
For $j\le M$ we use the Bernstein inequality
$$
\|\Delta_j f\|_{L^{2r}(\mathbb R^n)}^2 \le C \lambda_j^{\frac{n}{r}}\|\Delta_j f\|^2_{L^r(\mathbb R^n)}\le C \lambda_j^{\frac{n}{r}}\|f\|^2_{L^r(\mathbb R^n)}
$$
and so
$$
\ba
\sum_{-\infty}^{M}\|\Delta_j f\|_{L^{2r}(\mathbb R^n)}^2 \le C\|f\|_{L^r(\mathbb R^n)}^2\sum_{j=-\infty}^M\lambda_j^{\frac{n}{r}}\le \\
C\lambda_M^{\frac{n}{r}}\|f\|^2_{L^r(\mathbb R^n)}.
\ea
$$
Optimizing in $M$, we deduce (\ref{fro}).\\
\vspace{1cm}

{\bf{Acknowledment}} P.C.'s research was partially sponsored by NSF grant DMS-0804380. G. S.'s research was partially supported by  the
RFFI grant 08-01-00372-a. P.C. gratefully acknowledges the hospitality of Oxford University's OXPDE Center.

%G. Seregin\\
%Steklov Institute of Mathematics at St.Petersburg, \\
%St.Petersburg, Russia
%\\
%\\


\begin{thebibliography}{99}

%\bibitem {CKN}
%Caffarelli, L., Kohn, R.-V., Nirenberg, L., Partial regularity of
%suitable weak solutions of the Navier-Stokes equations, Comm. Pure
%Appl. Math., Vol. XXXV (1982), pp. 771--831.
%\bibitem {Cam}
%Campanato, S., Equazioni paraboliche del secondo ordine e spazi ${\mathcal L}^{2,\theta}(\Om,\delta)$. Ann. Mat. Pura
%Appl. 73(1966), 55-102.
\bibitem{bl} J. Bergh, J.  L\"{o}fstrom, Interpolation spaces, an introduction,
Springer-Verlag, Berlin (1976).
\bibitem{c} P. Constantin, Nonlinear Fokker-Planck Navier-Stokes Systems,
Commun. Math. Sci. {\bf 3} (2005), 531-544.

\bibitem{c-smo} P. Constantin, Smoluchowski Navier-Stokes systems, Contemporary Mathematics {\bf 429} G-Q Chen, E. Hsu, M. Pinsky editors, AMS, Providence (2007), 85-109.

\bibitem{cftz}  P. Constantin, C. Fefferman, E. Titi, A. Zarnescu, Regularity for coupled two-dimensional nonlinear Fokker-Planck and Navier-Stokes systems, Commun. Math. Phys., {\bf 270} (2007) 789-811.

\bibitem{cm} P. Constantin, N. Masmoudi, Global well-posedness
for a Smoluchowski equation coupled with Navier-Stokes equations in 2D, Commun. Math. Phys. {\bf 278} (2008), 179-191.

\bibitem{cs} P. Constantin, G. Seregin, H\"{o}lder Continuity of Solutions of 2D Navier-Stokes Equations with Singular Forcing, preprint.

%\bibitem{DaPr}
%Da Prato, G. Spazi ${\mathcal L}^{2,\theta}(\Om,\delta)$ e loro proprieta. Ann. Mat. Pura Appl. 69(1965), 383-392.

%\bibitem {GiaStr}
% Giaquinta, M., Struwe, M., On the Partial Regularity of Weak Solutions
%of Nonlinear Parabolic Systems, Math. Z. 179(1982), 437-451.

\bibitem{doied} M. Doi, S.F. Edwards, The Theory of Polymer Dynamics, Oxford
University Press, Oxford 1988.


\bibitem{La}
Ladyzhenskaya, O. A., Global solvability of a boundary value
roblem for the Navier-Stokes equations in the case of two spatial
variables, Doklady of the USSR, 123(1958), 427--429.

\bibitem{lemarie} P.G. Lemari\'{e}-Rieusset, Recent developments in the Navier-Stokes problem, Chapmann and Hall/CRC, Research Notes in Mathematics 431, CRC, Boca Raton, 2002.

\bibitem{lebris} C. LeBris, P - L. Lions, Existence and uniqueness of solutions to Fokker-Planck type equations with irregular coefficients, Commun. PDE {\bf 33}(7)(2008), 1272-1317.

\bibitem{fanghua} F. Lin, P. Zhang, Z. Zhang, On the global existence if smooth solution to the 2D FENE dumbell model, Commun. Math. Phys. {\bf 277} (2008), 531-553.




\end{thebibliography}
\end{document}